# An iterative algorithm for a common fixed point of Bregman Relatively Nonexpansive Mappings

Oladipo Abiodun Timoye[1] and Enyinnaya Ekuma-Okereke[2*]

[1]Department of Pure and Applied Mathematics, Ladoke Akintola University of Technology, Ogbomoso, P.M.B 4000, Nigeria. Email: atlab_3@yahoo.com
[2]Department of Mathematics/Computer Science, Federal University of Petroleum Resources, Effurun, Nigeria.
[*]Correspondence: ekuma.okereke@fupre.edu.ng

**Abstract**: We introduce and investigate an iterative scheme for approximating common fixed point of a family of Bregman relatively-nonexpansive mappings in real reflexive Banach spaces. We prove strong convergence theorem of the sequence generated by our scheme under some appropriate conditions. Furthermore, our scheme and results unify some known results obtained in this direction.

## 1. Introduction

Let $X$ denote a real reflexive Banach Space with the norm $\|.\|$, $X^*$ stands for the dual space of $X$. The normalized duality mapping from $X$ to $2^{X^*}$ denoted by $J$ is defined by

$$Jx = \{x^* \in X^*: \langle x, x^* \rangle = \|x\|^2 = \|x^*\|^2\}, \forall x \in X, \qquad (1.1)$$

where $\langle .,. \rangle$ denotes the generalized duality pairing between $X$ and $X^*$.

Let $C$ be a nonempty, closed and convex subsets of $X$. Let $T: C \to C$ be a nonlinear self mapping. $T$ is said to be nonexpansive mapping if $\|Tx - Ty\| \leq \|x - y\|, \forall x, y \in C$, and $T$ is said to be quasi-nonexpansive mapping if $\|Tx - p\| \leq \|x - p\|, \forall x \in C, p \in F(T)$, where $F(T) = \{x \in C: Tx = x\}$ is the set of fixed points of a mapping $T$. A point $p \in C$, is called an asymptotic fixed point of a mapping $T$ if $C$ contains a sequence $\{x_n\}$ with $x_n \rightharpoonup p$ such that $\|x_n - Tx_n\| = 0$. The set of asymptotic fixed point is denoted by $\hat{F}(T)$, (see [1]).

A mapping $T: C \to C$ is said to be

Bregman firmly nonexpansive (BFNE) (see [2]) if

$$\langle \nabla f(Tx) - \nabla f(Ty), Tx - Ty \rangle \leq \langle f(x) - \nabla f(y), Tx - Ty \rangle \forall x, y \in C,$$

or equivalently,

$$D_f(Tx,Ty) + D_f(Ty,Tx) + D_f(Tx,x) + D_f(Ty,y) \leq D_f(Tx,y) + D_f(Ty,x)$$

Bregman quasi-nonexpansive (BQNE) (see [3]) if $F(T) \neq \emptyset$ and

$$D_f(p,Tx) \leq D_f(p,x), \forall x \in C, \forall p \in F(T)$$

Bregman relatively-nonexpansive (BRNE) (see [3]) if $F(T) \neq \emptyset$ and

$$D_f(p,Tx) \leq D_f(p,x), \forall x \in C, \forall p \in F(T) = \hat{F}(T).$$

Existence and approximation of fixed points of relatively nonexpansive and quasi-nonexpansive mappings have extensively been studied by many authors for some decades now in Hilbert spaces, see for example, [4,5]. Since some of the methods fails to give same conclusion in Banach spaces which is more general than Hilbert spaces, for instance the resolvent $R_A = (1+A)^{-1}$ of a maximal monotone mapping $A: H \to 2^H$ and the metric projection $P_C$ onto a nonempty, closed and convex subset $C$ of $H$ are nonexpansive in Hilbert spaces but not nonexpansive in general Banach spaces. In this connection, Alber [6] introduced a generalized projection operator in Banach spaces which is an analogue of the metric projection in Hilbert spaces. Another way to overcome this problem is the use of distance function $D_f(.,.)$ introduced by Bregman [7] instead of norm now being studied by many authors. This has over the past seven years opened a growing area of research; see [8, 9].

Recently, in 2016, Alghamdi et al [20] introduced an iterative scheme for finding a common point of the fixed point se of a Bregman relatively nonexpansive mapping and the solution set of variational inequality problem for a continuous monotone mapping. They proved a strong convergence theorem for the sequences produced by the method.

In [19], Ugwunnadi and Ali proved a new strong convergence theorem for finite family of quasi-Bregman nonexpansive mapping and system of equilibrium problem in real Banach space.

In [12], Alghamdi et al proved a strong convergence theorem for the common fixed point of finite family of quasi-Bregman nonexpansive mappings.

Inspired and motivated by the researches ongoing in this direction, we consider an iterative scheme which converges strongly to a common fixed point of a family of Bregman relatively nonexpansive mappings in reflexive Banach spaces.

## 2. Preliminaries

Let $f: X \to (-\infty, +\infty]$ be assumed to be proper, lower-semicontinuous and convex function. Let the domain of $f$ be denoted as $dom f = \{x \in X: f(x) < +\infty\}$. Let $x \in int\, dom\, f$. The subdifferential of $f$ at $x$ is the convex set defined by

$$\partial f(x) = \{x^* \in X^*: f(x) + \langle x^*, y - x \rangle \leq f(y); \; \forall y \in X\}. \tag{2.1}$$

A function $f^*: X^* \to (-\infty, +\infty]$ defined by $f^*(x^*) = sup\{\langle x, x^* \rangle - f(x), \; x \in X\}$ is called the conjugate function of $f^*$. We see from the conjugate inequality that $f(x) \geq \langle x, x^* \rangle - f^*(x^*), \forall x \in X, \forall x^* \in X^*$, (see [14]). The function $f$ is said to be cofinite if $dom\, f^* = X^*$. A function $f$ on $X$ is coercive [11], if the sublevel set of $f$ is bounded, equivalently $\lim_{\|x\| \to \infty} f(x) = +\infty$. It is said to be strongly coercive [14], if $\lim_{\|x\| \to \infty} \frac{f(x)}{\|x\|} = +\infty$.

For any $x \in int\, dom\, f$ and $y \in X$, the right hand derivative of $f$ at $x$ in the direction of $y$ is defined by $f^\circ(x, y) = \lim_{t \to 0^+} \frac{f(x+ty) - f(x)}{t}$. A function $f$ is said to be Gâteaux differentiable at $x$ if $\lim_{t \to 0^+} \frac{f(x+ty) - f(x)}{t}$ exists for any $y$. In this case, $f^\circ(x, y)$ coincides with $\nabla f(x)$, the value of the gradient $\nabla f$ of $f$ at $x$. The function $f$ is said to be Gâteaux differentiable if it is Gâteaux differentiable for any $x \in int\, dom\, f$. The function $f$ is said to be Fréchet differentiable at $x$ if this limit is attained uniformly in $\|y\| = 1$. Finally, $f$ is said to be uniformly Fréchet differentiable on a subset $C$ of $X$ if the limit is attained uniformly for $x \in C$ and $\|y\| = 1$.

**Definition 2.1**: (Cf [15]), Let $f: X \to (-\infty, +\infty]$ be a Gâteaux differentiable function. The function $D_f: dom\, f \times int\, dom\, f \to [0, +\infty)$ defined by

$$D_f(y, x) = f(y) - f(x) - \langle \nabla f(x), y - x \rangle \tag{2.2}$$

is called the Bregman distance with respect to $f$. The Bregman distance has two important properties as follows:

Let $P_C^f: int\, dom\, f \to C$ be a mapping such that $P_C^f(x) \in C$ satisfying

$$D_f(P_C^f(x), x) = inf\{D_f(y, x): y \in C\} \tag{2.3}$$

is the Bregman Projection of $x \in int\, domf$ onto a nonempty closed and convex set $C \subset dom\, f$.

**Remark 2.2**: If $X$ is a smooth and strictly convex Banach spaces and $f(x) = \|x\|^2$ for all $x \in X$, then we have that $\nabla f(x) = 2Jx$, for all $x \in X$, where $J$ is the normalized duality mapping. Clearly, we obtain that

$$D_f(y, x) = f(y) - f(x) - \langle \nabla f(x), y - x \rangle$$
$$= \|y\|^2 - \|x\|^2 - 2\langle y, Jx \rangle + 2\|x\|^2$$

$$= \|x\|^2 - 2\langle y, Jx\rangle + \|y\|^2 = \phi(y, x) \quad \forall x, y \in X.$$

Which is Lyapunov function introduced by Alber [6] and $P_C^f(x)$ reduces to the generalized projection given as

$$\Pi_C(x) = \arg\min_{y \in C} \phi(y, x).$$

In addition, if $X$ coincides with $H$, in Hilbert space then $J = I$ and

$$D_f(y, x) = f(y) - f(x) - \langle \nabla f(x), y - x\rangle$$
$$= \|x\|^2 - \|y\|^2 - 2\langle x, y\rangle + 2\|y\|^2$$
$$= \|x\|^2 - 2\langle x, y\rangle + \|y\|^2 = \|x - y\|^2 \quad \forall x, y \in X.$$

Hence the Bregman Projection $P_C^f(x)$ reduces to metric projection of $H$ onto $C$, $P_C(x)$.

**Definition 2.3**: (see [18]), $f: X \to (-\infty, +\infty]$ is said to be Legendre function if it satisfies the following two conditions:

(L1) $int\ dom\ f \neq \emptyset, f$ is Gâteaux differentiable on $int\ dom\ f$ and $dom\ f = int\ dom\ f$,

(L2) $int\ dom\ f^* \neq \emptyset, f^*$ is Gâteaux differentiable on $int\ dom\ f^*$ and $dom\ f^* = int\ dom\ f^*$.

Remark 2.2: (cf [16-18]), since $X$ is reflexive, then we have that $(\partial f)^{-1} = \partial f^*$ and since $f$ is Legendre, then $\partial f$ is a bijection which satisfies $\nabla f = (\nabla f^*)^{-1}$, $ran\ \nabla f = dom\ \nabla f^* = int\ dom\ f^*$ and $ran\ \nabla f^* = dom\ \nabla f = int\ dom f$. $f$ and $f^*$ are strictly convex on their $int\ dom\ f$. If the subdifferential of $f$ is single valued, it coincides with the gradient of $f$, that is $\partial f = \nabla f$.

Example of a Legendre function is $f(x) = \frac{1}{p}\|x\|^p (1 < p < \infty)$. If $X$ is smooth and strictly convex Banach Spaces, then in this case the gradient $\nabla f$ coincides with the generalized duality mapping of $X$, that is $\nabla f = J_p$. If the space is a Hilbert space, $H$, then $\nabla f = I$, where $I$ is the identity mapping in $H$. Throughout this paper, we assumed that $f$ is Legendre.

**Definition 2.4** Let $f: X \to (-\infty, +\infty]$ be a $G\dot{a}teaux$ differentiable function. The modulus of total convexity of $f$ at $x \in intdom\ f$ is the function $V_f(x, .): intdom\ f \times [0, +\infty) \to [0, +\infty)$ defined by

$$V_f(x, t) = \inf\{D_f(y, x): y \in dom\ f, \|y - x\| = t\}. \tag{2.4}$$

The function $f$ is called totally convex at $x$ if $V_f(x, t) > 0$ whenever $t > 0$. The function $f$ is called totally convex if it is totally convex at any point $x \in int\ dom\ f$. The function is said to be totally convex on bounded sets if $V_f(B, t) > 0$ for any nonempty bounded subset $B$ of $X$ and $t > 0$, where the modulus of totall convexity of the function $f$ on the set $B$ is the function $V_f: int\ domf \times [0, +\infty) \to [0, +\infty)$ defined by

$$V_f(B, t) = \inf\{V_f(x, t): x \in B \cap dom\ f\}. \tag{2.5}$$

Also in this paper, we shall make use of the function $V_f: X^* \times X \to [0, +\infty)$ associated with $f$ defined by

$$V_f(x^*, x) = f(x^*) - \langle x^*, x \rangle + f^*(x), \forall x \in X, x^* \in X^*. \tag{2.6}$$

We see that $V_f(,) \geq 0$ and the relation

$$V_f(x^*, x) = D_f(\nabla f^*(x^*), x) \tag{2.7}$$

Moreover, by the subdifferential inequality, we obtain

$$V_f(x^*, x) + \langle y^*, \nabla f^*(x^*) - x \rangle \leq V_f(x^* + y^*, x), \ \forall x \in X\ and\ x^*, y^* \in X^*. \tag{2.8}$$

We remark that $V_f$ is convex in the first variable.

In the sequel, we shall make use of the following lemmas

**Lemma 2.5** (see [22]). The function $f$ is totally convex on bounded sets if and only if for any two sequences $\{x_n\}$ and $\{y_n\}$ in $X$ such that the first one is bounded, then

$$\lim_{n \to \infty} D_f(y_n, x_n) = 0 \Longrightarrow \|y_n - x_n\| = 0.$$

**Lemma 2.6** (see [2]). Let $C$ be a nonempty, closed and convex subsets of $int\ dom\ f$ and $T: C \to C$ be a quasi-Bregman nonexpansive mapping with respect to $f$. Then $F(T)$ is closed and convex.

**Lemma 2.7** (see [23]). Let $C$ be a nonempty, closed and convex subsets of $X$. Let $f: X \to (-\infty, +\infty]$ $G\hat{a}teaux$ differentiable and totally convex function and let $x \in X$, then

$$z = P_C^f(x)\ if\ and\ if\ \langle \nabla f(x) - \nabla f(z), y - z \rangle \leq 0, \forall y \in C.$$
$$D_f\big(y, P_C^f(x)\big) + D_f\big(P_C^f(x), x\big) \leq D_f(y, x)\ \forall y \in C.$$

**Lemma 2.8** (see [11]). Let $X$ be a reflexive Banach space and let $f: X \to R$ be a continuous convex function which is strongly coercive. Then the following assertions are equivalent:

(1) $f$ is bounded on bounded subsets and uniformly smooth on bounded subsets of $X$.

(2) $f^*$ is Fréchet differentiable and $f^*$ is uniformly norm-to-norm continuous on bounded subsets of $X^*$.

(3) $domf^* = X^*, f^*$ is strongly coercive and uniformly convex on bounded subsets of $X^*$.

**Lemma 2.9** (see [13]). Let $X$ be a Banach space, let $r > 0$ be a constant and let $f: X \to R$ be a continuous and convex function which is uniformly convex on bounded subsets of $X$. Then

$$f(\sum_{k=0}^{\infty} \alpha_k x_k) \leq \sum_{k=0}^{\infty} \alpha_k f(x_k) - \alpha_i \alpha_j \rho_r(\|x_i - x_j\|),$$

for all $i, j \in \mathbb{N} \cup \{0\}$, $x_k \in B_r$, $\alpha_k \in (0,1)$ and $k \in \mathbb{N} \cup \{0\}$ with $\sum_{k=0}^{\infty} \alpha_k = 1$, where $\rho_r$ is the gauge of uniform convexity of $f$.

**Lemma 2.10** (see [21]). If $f: X \to (-\infty, +\infty]$ is uniformly Fréchet differentiable and bounded on bounded subsets of $X$, then $\nabla f$ is uniformly continuous on bounded subsets of $X$ from the strong topology of $X$ to the strong topology of $X^*$.

**Lemma 2.11** (see [10]) Let $f: X \to (-\infty, +\infty]$ be a Gâteaux differentiable and totally convex function if $x_0 \in X$ and the sequence $\{D_f(x_n, x_0)\}$ is bounded, then the sequence $\{x_n\}$ is also bounded.

**Lemma 2.12** (see [17]). Let $f: X \to (-\infty, +\infty]$ be a proper, lower semi-continuous and convex function, then $f^*: X^* \to (-\infty, +\infty]$ is a proper, weak* lower semi-continuous and convex function. Thus, for all $z \in X$, we have

$$D_f(z, \nabla f^*(\sum_{i=1}^{N} t_i \nabla f(x_i))) \leq \sum_{i=1}^{N} t_i D_f(z, x_i)$$

**Lemma 2.13** (see [24]). Let $\{a_n\}_{n=1}^{\infty}$ be a sequence of nonnegative real numbers satisfying the following relation: $a_{n+1} \leq (1 - \alpha_n) a_n + \alpha_n \delta_n$, $n \geq n_0$, where $\{\alpha_n\}_{n=1}^{\infty}$ is a sequence in $(0,1)$, $\{\delta_n\}_{n=1}^{\infty}$ is a sequence in $R$ satisfying the following conditions: $\lim_{n \to \infty} \alpha_n = 0$, $\sum_{n=1}^{\infty} \alpha_n = \infty$, $\lim_{n \to \infty} \sup \delta_n \leq 0$. Then $\lim_{n \to \infty} a_n = 0$.

**Lemma 2.14** (see [25]). Let $\{a_n\}_{n=1}^{\infty}$ be a sequence of nonnegative real numbers such that there exists a nondecreasing subsequence $\{n_i\}$ of $\{n\}$ that is $a_{n_i} \leq a_{n_{i+1}}$ for all $i \in N$. Then there exists a nondecreasing subsequence $\{m_k\} \subset N$ such that $m_k \to \infty$ and the following properties are satisfied for all (sufficiently large number $k \in N$): $a_{m_k} \leq a_{m_{k+1}}$ and $a_k \leq a_{m_{k+1}}$. In fact, $m_k = max\{j \leq k : a_j \leq a_{j+1}\}$.

## 3. Main Results

**Theorem 3.1**: Let $C$ be a nonempty, closed and convex subset of $intdom\ f$, let $f : X \to R$ be a strongly coercive Legendre function which is bounded, uniformly Fréchet differentiable and totally convex on bounded subsets of a real reflexive Banach Space $X$. Let $T_1, T_2: C \to C$ be a family of Bregman relatively nonexpansive mappings. Assume that $\mathcal{F} = F(T_1) \cap F(T_2) \neq \emptyset$.

For any fixed $u, x_0 \in C$, let $\{x_n\}$ be a sequence of $C$ generated by the following iterative algorithm:

$$\begin{cases} z_n = \nabla f^*(c_n \nabla f(x_n) + (1-c_n)\nabla f(T_2 x_n)); \\ y_n = \nabla f^*(\beta_n \nabla f(x_n) + (1-\beta_n)\nabla f(T_1 x_n)); \\ x_{n+1} = P_C^f \nabla f^*(\alpha_n \nabla f(u) + (1-\alpha_n)(\theta_n \nabla f(x_n) + \delta_n \nabla f(y_n) + \gamma_n \nabla f(z_n))), \quad n \geq 0, \end{cases} \quad (3.1)$$

where $\{\delta_n\}, \{\theta_n\}, \{\gamma_n\}$ are sequences in $(0,1)$, $\{\alpha_n\}$ is a sequence in $(0,1)$ satisfying the following conditions: $(i) \lim_{n\to\infty} \alpha_n = 0$ $(ii) \sum_{n=1}^{\infty} \alpha_n = \infty$ $(iii) \theta_n + \delta_n + \gamma_n = 1$. Then, $\{x_n\}$ converges strongly to a common fixed point of $T_1$ and $T_2$ nearest to $u$.

*Proof.* Now by Lemma 2.6, we obtain that $\mathcal{F}$ is closed and convex. Let $p = P_{\mathcal{F}}^f \in \mathcal{F}$.
Now setting
$w_n = \nabla f^*(\theta_n \nabla f(x_n) + \delta_n \nabla f(y_n) + \gamma_n \nabla f(z_n))$ and
$h_n = \nabla f^*(\alpha_n \nabla f(u) + (1-\alpha_n)\nabla f(w_n))$, then $x_{n+1} = P_C^f h_n$. (3.2)

Now, from Lemma 2.8 and since $f$ is bounded and uniformly smooth on bounded subsets of $X$, so $f^*$ is uniformly convex on bounded subsets of $X^*$. Then using Lemma 2.9, the properties of $D_f$ and $T_1, T_2$, and from (3.1), (2.6), (2.7) we obtain that

$$\begin{aligned} D_f(p, y_n) &= D_f\left(p, \nabla f^*(\beta_n \nabla f(x_n) + (1-\beta_n)\nabla f(T_1 x_n))\right) \\ &= V_f(p, \beta_n \nabla f(x_n) + (1-\beta_n)\nabla f(T_1 x_n)) \\ &\leq f(p) - \langle p, \beta_n \nabla f(x_n) + (1-\beta_n)\nabla f(T_1 x_n) \rangle \\ &\quad + f^*(\beta_n \nabla f(x_n) + (1-\beta_n)\nabla f(T_1 x_n)) \\ &\leq \beta_n f(p) + (1-\beta_n)f(p) - \beta_n \langle p, \nabla f(x_n) \rangle + (1-\beta_n)\langle p, \nabla f(T_1 x_n) \rangle \\ &\quad + \beta_n f^*(\nabla f(x_n)) + (1-\beta_n)f^*(\nabla f(T_1 x_n)) \\ &\quad - \beta_n(1-\beta_n) p_s^*(\|\nabla f(x_n) - \nabla f(T_1 x_n)\|) \\ &= \beta_n V_f(p, \nabla f(x_n)) + (1-\beta_n)V_f(p, \nabla f(T_1 x_n)) \\ &\quad - \beta_n(1-\beta_n) p_s^*(\|\nabla f(x_n) - \nabla f(T_1 x_n)\|) \\ &= \beta_n D_f(p, x_n) + (1-\beta_n)D_f(p, T_1 x_n) \\ &\quad - \beta_n(1-\beta_n) p_s^*(\|\nabla f(x_n) - \nabla f(T_1 x_n)\|) \\ D_f(p, y_n) &\leq D_f(p, x_n) - \beta_n(1-\beta_n) p_s^*(\|\nabla f(x_n) - \nabla f(T_1 x_n)\|) \end{aligned} \quad (3.3)$$

$$\leq D_f(p, x_n).$$

Similarly,
$$D_f(p, z_n) \leq D_f(p, x_n) - c_n(1-c_n)p_s^*(\|\nabla f(x_n) - \nabla f(T_2 x_n)\|) \tag{3.4}$$
$$\leq D_f(p, x_n).$$

In addition, employing (2.6), (2.7), we obtain

$$\begin{aligned} D_f(p, w_n) &= D_f\left(p, \nabla f^*(\theta_n \nabla f(x_n) + \delta_n \nabla f(y_n) + \gamma_n \nabla f(z_n))\right) \\ &= V_f(p, \theta_n \nabla f(x_n) + \delta_n \nabla f(y_n) + \gamma_n \nabla f(z_n)) \\ &\leq f(p) - \langle p, \theta_n \nabla f(x_n) + \delta_n \nabla f(y_n) + \gamma_n \nabla f(z_n) \rangle \\ &\quad + f^*(\theta_n \nabla f(x_n) + \delta_n \nabla f(y_n) + \gamma_n \nabla f(z_n)) \\ &\leq f(p) - \theta_n \langle p, \nabla f(x_n) \rangle - \delta_n \langle p, \nabla f(y_n) \rangle - \gamma_n \langle p, \nabla f(z_n) \rangle \\ &\quad + \theta_n f^*(\nabla f(x_n)) + \delta_n f^*(\nabla f(y_n)) + \gamma_n f^*(\nabla f(z_n)) \\ &= \theta_n \left(f(p) - \langle p, \nabla f(x_n) \rangle + f^*(\nabla f(x_n))\right) \\ &\quad + \delta_n \left(f(p) - \langle p, \nabla f(y_n) \rangle + f^*(\nabla f(y_n))\right) \\ &\quad + \gamma_n \left(f(p) - \langle p, \nabla f(z_n) \rangle + f^*(\nabla f(z_n))\right) \\ &= \theta_n V_f(p, \nabla f(x_n)) + \delta_n V_f(p, \nabla f(y_n)) + \gamma_n V_f(p, \nabla f(z_n)) \\ &= \theta_n D_f(p, x_n) + \delta_n D_f(p, y_n) + \gamma_n D_f(p, z_n). \end{aligned} \tag{3.5}$$

Substituting (3.3) and (3.4) into (3.5), we obtain
$$\begin{aligned} D_f(p, w_n) &\leq D_f(p, x_n) - \delta_n \beta_n(1-\beta_n)p_s^*(\|\nabla f(x_n) - \nabla f(T_1 x_n)\|) \\ &\quad - \gamma_n c_n(1-c_n)p_s^*(\|\nabla f(x_n) - \nabla f(T_2 x_n)\|). \end{aligned} \tag{3.6}$$

Furthermore, from Lemma 2.7, (3.2), (3.6) and the property of $D_f$, we obtain
$$\begin{aligned} D_f(p, x_{n+1}) &= D_f\left(p, P_C^f h_n\right) \\ &\leq D_f(p, h_n) \\ &= D_f\left(p, \nabla f^*(\alpha_n \nabla f(u) + (1-\alpha_n)\nabla f(w_n))\right) \\ &\leq \alpha_n D_f(p, u) + (1-\alpha_n)D_f(p, w_n) \\ &\leq \alpha_n D_f(p, u) + (1-\alpha_n)D_f(p, x_n) \\ &\quad - (1-\alpha_n)\delta_n \beta_n(1-\beta_n)p_s^*(\|\nabla f(x_n) - \nabla f(T_1 x_n)\|) \\ &\quad - (1-\alpha_n)\gamma_n c_n(1-c_n)p_s^*(\|\nabla f(x_n) - \nabla f(T_2 x_n)\|) \tag{3.7} \\ &\leq \alpha_n D_f(p, u) + (1-\alpha_n)D_f(p, x_n). \tag{3.8} \end{aligned}$$

Thus by induction, we obtain that

$$D_f(p, x_{n+1}) \leq \max\{D_f(p, u), D_f(p, x_0)\}, \forall n \geq 0,$$

which implies that $\{D_f(p, x_n)\}$ and hence $\{D_f(p, T_1 x_n)\}$ are bounded. Thus we get from Lemmas 2.10, 2.11 that $\{x_n\}, \{y_n\}, \{z_n\}, \{w_n\}$ and $\{h_n\}$ are all bounded.

Furthermore, from (3.2), Lemma 2.7, (2.7) and (2.8), we obtain

$$\begin{aligned}
D_f(p, x_{n+1}) &= D_f(p, P_C^f h_n) \\
&\leq D_f(p, h_n) \\
&= D_f\left(p, \nabla f^*(\alpha_n \nabla f(u) + (1-\alpha_n)\nabla f(w_n))\right) \\
&= V_f(p, \alpha_n \nabla f(u) + (1-\alpha_n)\nabla f(w_n)) \\
&\leq V_f\left(p, \alpha_n \nabla f(u) + (1-\alpha_n)\nabla f(w_n) - \alpha_n(\nabla f(u) - \nabla f(p))\right) \\
&\quad + \alpha_n \langle \nabla f(u) - \nabla f(p), h_n - p \rangle \\
&= V_f(p, \alpha_n \nabla f(p) + (1-\alpha_n)\nabla f(w_n)) + \alpha_n \langle \nabla f(u) - \nabla f(p), h_n - p \rangle \\
&= D_f\left(p, \nabla f^*(\alpha_n \nabla f(p) + (1-\alpha_n)\nabla f(w_n))\right) + \alpha_n \langle \nabla f(u) - \nabla f(p), h_n - p \rangle \\
&\leq \alpha_n D_f(p, p) + (1-\alpha_n) D_f(p, w_n) + \alpha_n \langle \nabla f(u) - \nabla f(p), h_n - p \rangle \\
&\leq (1-\alpha_n) D_f(p, w_n) + \alpha_n \langle \nabla f(u) - \nabla f(p), h_n - p \rangle
\end{aligned}$$

$$\begin{aligned}
D_f(p, x_{n+1}) &\leq (1-\alpha_n) D_f(p, x_n) - (1-\alpha_n)\delta_n \beta_n (1-\beta_n) p_s^*(\|\nabla f(x_n) - \nabla f(T_1 x_n)\|) \\
&\quad -(1-\alpha_n)\gamma_n c_n (1-c_n) p_s^*(\|\nabla f(x_n) - \nabla f(T_2 x_n)\|) \\
&\quad + \alpha_n \langle \nabla f(u) - \nabla f(p), h_n - p \rangle \quad (3.9) \\
&\leq (1-\alpha_n) D_f(p, x_n) + \alpha_n \langle \nabla f(u) - \nabla f(p), h_n - p \rangle. \quad (3.10)
\end{aligned}$$

We now consider two cases.

**Cases I:** Suppose that there exists $n_0 \in N$ such that $\{D_f(p, x_n)\}$ is monotone non-increasing for all $n \geq n_0$. Then we get that $\{D_f(p, x_n)\}$ is convergent and $D_f(p, x_n) - D_f(p, x_{n+1}) \to 0$, so that from (3.9), we obtain that

$$(1-\alpha_n)\delta_n \beta_n (1-\beta_n) p_s^*(\|\nabla f(x_n) - \nabla f(T_1 x_n)\|) \to 0, \quad (3.11)$$

and

$$(1-\alpha_n)\gamma_n c_n (1-c_n) p_s^*(\|\nabla f(x_n) - \nabla f(T_2 x_n)\|) \to 0 \quad (3.12)$$

which by the property of $p_s^*$ give

$$\nabla f(x_n) - \nabla f(T_1 x_n) \to 0, \quad \nabla f(x_n) - \nabla f(T_2 x_n) \to 0 \text{ as } n \to \infty. \quad (3.13)$$

Moreover, from (3.1) and (3.13), we obtain

$$\|\nabla f(y_n) - \nabla f(x_n)\| \leq \beta_n \|\nabla f(x_n) - \nabla f(x_n)\| - (1-\beta_n)\|\nabla f(T_1 x_n) - \nabla f(x_n)\| \to 0, \quad (3.14)$$

$$\|\nabla f(z_n) - \nabla f(x_n)\| \leq c_n \|\nabla f(x_n) - \nabla f(x_n)\| - (1-c_n)\|\nabla f(T_2 x_n) - \nabla f(x_n)\| \to 0. \quad (3.15)$$

In addition, employing (3.1) and (3.14), (3.15) and (3.13) we obtain

$$\|\nabla f(w_n) - \nabla f(x_n)\| \leq \theta_n \|\nabla f(x_n) - \nabla f(x_n)\| + \delta_n \|\nabla f(y_n) - \nabla f(x_n)\|$$
$$+ \gamma_n \|\nabla f(z_n) - \nabla f(x_n)\|$$
$$\leq \delta_n \beta_n \|\nabla f(x_n) - \nabla f(x_n)\| - \delta_n (1-\beta_n)\|\nabla f(T_1 x_n) - \nabla f(x_n)\| \quad (3.16)$$
$$+ \gamma_n c_n \|\nabla f(x_n) - \nabla f(x_n)\| - \gamma_n (1-c_n)\|\nabla f(T_2 x_n) - \nabla f(x_n)\| \to 0.$$

Since, $f$ is strongly coercive and uniformly convex on bounded subsets of $X$, $f^*$ is uniformly Fréchet differentiable on bounded subsets of $X^*$ and by Lemma 2.8 we get that $\nabla f^*$ is uniformly continuous. So by this and together with (3.13), (3.14) and (3.15), we obtain that

$$x_n - T_1 x_n \to 0, \; x_n - T_2 x_n \to 0, \; x_n \to y_n \to 0, \; x_n \to z_n \to 0, \; x_n - w_n \to 0 \text{ as } n \to \infty. \quad (3.17)$$

Moreso, from Lemma 2.12 and (i), we obtain that

$$D_f(w_n, h_n) = D_f\left(w_n, \nabla f^*(\alpha_n \nabla f(u) + (1-\alpha_n)\nabla f(w_n))\right)$$
$$\leq \alpha_n D_f(w_n, u) + (1-\alpha_n) D_f(w_n, w_n) \to 0 \text{ as } n \to \infty, \quad (3.18)$$

and by Lemma 2.5, we obtain that

$$w_n - h_n \to 0 \text{ as } n \to \infty. \quad (3.19)$$

Now since $X$ is reflexive and $\{h_n\}$ is bounded, there exists a subsequence $\{h_{n_i}\}$ of $\{h_n\}$ such that $h_{n_i} \rightharpoonup h \in C$, and

$$\lim_{n\to\infty} sup \langle \nabla f(u) - \nabla f(p), h_n - p \rangle = \lim_{i\to\infty} sup \langle \nabla f(u) - \nabla f(p), h_{n_i} - p \rangle.$$

Hence, we obtain from (3.19) and (3.17), that $x_{n_i} \rightharpoonup h$. Using (3.17) and the fact that $T_1, T_2$ are Bregman relatively nonexpansive mappings, we obtain that $h \in F(T)$ and by Lemma 2.7

$$\lim_{n\to\infty} sup \langle \nabla f(u) - \nabla f(p), h_n - p \rangle = \lim_{i\to\infty} sup \langle \nabla f(u) - \nabla f(p), h_{n_i} - p \rangle$$
$$= \langle \nabla f(u) - \nabla f(p), h - p \rangle \leq 0. \quad (3.20)$$

It therefore follows from (3.10), (3.20) and Lemma 2.13, that $D_f(p, x_n) \to 0 \; as \; n \to \infty$. Consequently, from Lemma 2.5, we obtain that $x_n \to p = P_\mathcal{F}^f(u)$.

**Case II:** Suppose that there exists a subsequence $\{n_i\}$ of $\{n\}$ such that

$$D_f(p, x_{n_i}) < D_f(p, x_{n_{i+1}}) \text{ for all } i \in N. \quad (3.21)$$

Then by lemma 2,14, there exists a nondecreasing sequence $\{m_k\} \subset N$ such that $m_k \to \infty$ and

$D_f(p, x_{m_k}) \leq D_f(p, x_{m_{k+1}})$, $D_f(p, x_k) \leq D_f(p, x_{m_{k+1}})$ for all $k \in N$. Then from (3.9) and the fact that $\alpha_{m_k} \to 0$, we obtain that

$$p_s^*(\|\nabla f(x_{m_k}) - \nabla f(T_1 x_{m_k})\|) \to 0, \text{ and } (\|\nabla f(x_{m_k}) - \nabla f(T_2 x_{m_k})\|) \to 0 \text{ as } k \to \infty. \quad (3.22)$$

Thus we get from the same method of proof in **Case I** that

$$x_{m_k} - T_1 x_{m_k} \to 0, \quad x_{m_k} - T_2 x_{m_k} \to 0, x_{m_k} \to y_{m_k} \to 0, \; x_{m_k} \to z_{m_k} \to 0 \text{ as } k \to \infty, \quad (3.23)$$

$x_{m_k} - w_{m_k} \to 0$ as $k \to \infty$

and also we obtain

$$\lim_{k \to \infty} \sup \langle \nabla f(u) - \nabla f(p), h_{m_k} - p \rangle \leq 0 \quad (3.24)$$

Now from (3.10), we have that

$$D_f(p, x_{m_{k+1}}) \leq (1 - \alpha_{m_k}) D_f(p, x_{m_k}) + \alpha_{m_k} \langle \nabla f(u) - \nabla f(p), h_{m_k} - p \rangle. \quad (3.25)$$

Since $D_f(p, x_{m_k}) \leq D_f(p, x_{m_{k+1}})$, we have

$$\alpha_{m_k} D_f(p, x_{m_k}) \leq D_f(p, x_{m_k}) - D_f(p, x_{m_{k+1}}) + \alpha_{m_k} \langle \nabla f(u) - \nabla f(p), h_{m_k} - p \rangle.$$

Therefore,

$$\alpha_{m_k} D_f(p, x_{m_k}) \leq \alpha_{m_k} \langle \nabla f(u) - \nabla f(p), h_{m_k} - p \rangle. \quad (3.26)$$

Using (3.24), then (3.26) implies

$$D_f(p, x_{m_k}) \to 0 \text{ as } k \to \infty. \quad (3.27)$$

Consequently,

$$D_f(p, x_{m_{k+1}}) \to 0 \text{ as } k \to \infty. \quad (3.28)$$

But $D_f(p, x_k) \leq D_f(p, x_{m_{k+1}})$ for all $k \in N$. Thus we obtain that $D_f(p, x_k) \to 0$ as $n \to \infty$. Hence, by lemma 2.5, we have that $x_k \to p$. Therefore, from the above Cases, we can conclude that $\{x_n\}$ converges strongly to a common fixed point of $T_1$ and $T_2$ which is $p = P_{\mathcal{F}}^f(u)$ and that completes the proof of our theorem.

We observe that the method of proof of Theorem 3.1 provides a convergence theorem for a finite family of Bregman relatively nonexpansive mappings. The following theorem suffices.

**Theorem 3.2:** Let $C$ be a nonempty, closed and convex subset of $intdom \, f$, let $f : X \to R$ be a strongly coercive Legendre function which is bounded, uniformly Fréchet differentiable and totally convex on bounded subsets of a real reflexive Banach Space $X$. Let $T_i : C \to C, i = 1, 2, \ldots, N$ be a family of Bregman relatively nonexpansive mappings. Assume that $\mathcal{F} =$

$\cap_{i=1}^{N} F(T_i) \neq \emptyset$. For any fixed $u, x_0 \in C$, let $\{x_n\}$ be a sequence of $C$ generated by the following iterative algorithm:

$$\begin{cases} y_n = \nabla f^*(\beta_n \nabla f(x_n) + (1-\beta_n) \nabla f(T_i x_n)), & i = 1,2,\ldots,N; \\ x_{n+1} = P_C^f \nabla f^*\left(\alpha_n \nabla f(u) + (1-\alpha_n)\left(\theta_{n,0} \nabla f(x_n) + \sum_{i=1}^{N} \theta_{n,i} \nabla f(y_n)\right)\right), & n \geq 0 \end{cases} \quad (3.29)$$

where $\{\theta_{ni}: i = 0,1,2,\ldots,N\}$ are sequences in $(0,1)$, $\{\alpha_n\}$ is a sequence in $(0,1)$ satisfying the following conditions: $(i)\ \lim_{n\to\infty} \alpha_n = 0$ $(ii)\ \sum_{n=1}^{\infty} \alpha_n = \infty$ $(iii)\ \theta_{n,0} + \sum_{i=1}^{N} \theta_{n,i} = 1$.

Then, $\{x_n\}$ converges strongly to a common fixed point of $T_i, i = 1,2,\ldots,N$ nearest to $u$.

**Competing interests**

The authors declare that they have no competing interests.

**Acknowledgement**

The authors are grateful to the referees for their careful reading and suggestions.